# Bias correction and confidence intervals following sequential tests

Tze Leung Lai[1],[*]  Zheng Su[1] and Chin Shan Chuang[2]

*Stanford University and Millennium Partners*

**Abstract:** An important statistical inference problem in sequential analysis is the construction of confidence intervals following sequential tests, to which Michael Woodroofe has made fundamental contributions. This paper reviews Woodroofe's method and other approaches in the literature. In particular it shows how a bias-corrected pivot originally introduced by Woodroofe can be used as an improved root for sequential bootstrap confidence intervals.

## 1. Introduction and overview

Estimation following sequential tests is an important but difficult problem in sequential analysis. When the sample size is not fixed in advance but is a random variable $T$ that depends on the data collected so far, the sample moments and maximum likelihood estimates of population parameters can have substantial bias. For example, while the sample mean $\bar{X}_n$ is an unbiased estimate of the population mean $\mu$ based on a sample of $n$ i.i.d. observations $X_1, \ldots, X_n$, $\bar{X}_T$ is clearly biased upwards if $T$ is the first time when $S_n = n\bar{X}_n$ exceeds some threshold. How can one correct for the bias due to data-dependent sample size and how should one construct a confidence interval for $\mu$? Following Siegmund's seminal paper [18] on this problem, there have been many important developments in the literature, among which are the fundamental contributions of Woodroofe and his collaborators.

The simplest method to construct a confidence interval for the population mean $\mu$ is to use the naive normal approximation that treats $\sqrt{T}(\bar{X}_T - \mu)/\hat{\sigma}_T$ as approximately standard normal, where $\hat{\sigma}_n^2$ is a consistent estimate of $\text{Var}(X_1)$. Although one can justify the normal approximation by appealing to Anscombe's theorem [1] when

$$(1.1) \qquad T/a \xrightarrow{P} c \quad \text{as} \quad a \to \infty$$

for some nonrandom constant $c > 0$ and design parameter $a > 0$, the normal approximation essentially treats $T$ as nonrandom and has been found to be unsatisfactory in practice. This led Siegmund to develop exact methods for constructing confidence intervals for $\mu$ in the normal case by using a certain ordering of the sample space of $(T, S_T)$ when $T$ is the stopping time of a repeated significance test. Earlier Armitage used [3] numerical methods to evaluate exact confidence intervals for a Bernoulli parameter following sequential tests. Siegmund's approach was

---

[*]This research was supported by the National Science Foundation Grant DMS-0305749.

[1]Department of Statistics, Stanford Unversity, Stanford, CA 94305–4065, USA, e-mail: `lait@stat.stanford.edu`; `zhengsu@stanford.edu`

[2]Millennium Partners, e-mail: `cchuang@mlp.com`

*AMS 2000 subject classifications:* primary 62L12, 62G09, 62G20; secondary 60F05.

*Keywords and phrases:* bias due to optional stopping, bootstrap, coverage probabilities, hybrid resampling, very weak expansions.

subsequently extended to other stopping rules and to alternative orderings of the sample space by Tsiatis, Rosner and Mehta [20], Chang and O'Brien [5], Rosner and Tsiatis [17] and Emerson and Fleming [12].

Woodroofe [22] introduced "very weak" asymptotic expansions to correct for Anscombe's normal approximation in a one-parameter exponential family of densities $f_\theta(x) = e^{\theta x - \psi(\theta)}$ with natural parameter space $\Theta$. Denoting the stopping rule $T$ by $T_a$ to indicate its dependence on $a$, he strengthened (1.1) into

$$(1.2a) \qquad \lim_{a \to \infty} E_\theta |a/T_a - \kappa(\theta)| = 0 \text{ for a.e. } \theta \in \Theta^o,$$

$$(1.2b) \qquad \lim_{a \to \infty} a \int_C P_\theta\{T_a \leq a\eta_C\}\, d\theta = 0$$
$$\text{for some } \eta_C > 0 \text{ and every compact } C \subset \Theta^o,$$

where $\Theta^o$ denotes the interior of $\Theta$ and $\kappa : \Theta^o \to \mathbf{R}$ is continuous. Let $\mu = \dot\psi(\theta)$ and $\sigma^2 = \ddot\psi(\theta)$. Let $\widehat\theta_n$ be the maximum likelihood estimate of $\theta$ based on $X_1, \ldots, X_n$, and $\widehat\mu_n = \dot\psi(\widehat\theta_n) = \bar X_n$, $\widehat\sigma_n^2 = \ddot\psi(\widehat\theta_n)$. Consider a Bayesian model in which $\theta$ has a twice continuously differentiable prior density function $\xi$ with compact support $[\theta_0, \theta_1] \subset \Theta^o$. With $T = T_a$ satisfying (1.2a,b), Woodroofe [22] developed an asymptotic expansion for

$$(1.3) \qquad \int P_\theta\{\sqrt{T}(\bar X_T - \mu)/\widehat\sigma_T \leq c + b(\widehat\theta_T)/\sqrt{T}\}\xi(\theta)\, d\theta,$$

where $b$ is a piecewise continuous function on $\Theta^o$, and used it to construct a confidence interval $I$ for $\mu$ whose integrated coverage error

$$(1.4) \qquad \int P_\theta\{\mu \notin I\}\xi(\theta)\, d\theta$$

differs from the nominal value $2\alpha$ by $o(a^{-1})$. Subsequently, Woodroofe [23, 25, 26] showed how a version of Stein's identity [19] could be used to derive expressions for posterior expectations from which asymptotic expansions for (1.3) could be readily guessed. Moreover, for normal $X_i$ with known variance 1, Woodroofe [25] made use of these expansions to derive bias-corrected and renormalized pivots of the form

$$(1.5) \qquad R(\mu) = \{\sqrt{T}(\bar X_T - \mu) - T^{-1/2}b(\bar X_T)\}/\{1 + (2T)^{-1}b^2(\bar X_T)\},$$

where $b(\mu) = (\sqrt{\kappa(\mu)})'/\sqrt{\kappa(\mu)} = \kappa'(\mu)/\{2\kappa(\mu)\}$, in which $\kappa$ is given by (1.2a), noting that $\mu = \dot\psi(\theta) = \theta$ in the normal case.

Instead of using the Wald-type statistic $\sqrt{T}(\bar X_T - \mu)/\widehat\sigma_T$ as in (1.3), Coad and Woodroofe [8] and Weng and Woodroofe [21] considered confidence intervals based on signed-root likelihood ratio statistics, for which they developed very weak expansions leading to renormalized and bias-corrected signed-root likelihood ratio statistics as pivotal quantities. Woodroofe and his collaborators have also developed very weak asymptotic expansions to construct confidence sets in other sequential/adaptive experiments; see [9, 10, 23, 24, 27, 28].

The preceding methods assume parametric models, and more precisely, exponential families. For samples of fixed size, an important methodology for constructing confidence intervals without distributional assumptions is Efron's bootstrap method [11]. The bootstrap methodology can be extended as follows from the fixed sample size case to the case where the sample size is determined by a stopping rule $T$. Let $X_1, X_2, \ldots$ be i.i.d. random variables with a common unknown distribution function $F$. Given a randomly stopped sample $(X_1, \ldots, X_T)$, let $\widehat{F}$ denote the empirical distribution that puts probability mass $1/T$ at each of the sample values $X_i$. Let $X_1^*, X_2^*, \ldots$ be i.i.d. random variables with common distribution $\widehat{F}$ and let $T^*$ denote the corresponding stopping time for the sequence $\{X_i^*\}$. The *sequential bootstrap sample* $(X_1^*, \ldots, X_{T^*}^*)$ can be used to construct confidence intervals as in the nonsequential case, and Chuang and Lai [6] have studied the coverage probabilities of these bootstrap confidence intervals in the setting where $T$ is the stopping rule of a group sequential test, for which stopping can only occur at a prespecified subset $\{n_1, \ldots, n_K\}$ of positive integers, and have shown that the bootstrap method does not yield reliable confidence intervals because $\sqrt{T}(\bar{X}_T - \mu)/\widehat{\sigma}_T$ is no longer an approximate pivot. There are "pockets" of the parameter space where $T/a$ has a nondegenerate limiting distribution that varies with $\mu$, thus violating (1.1) and making the distribution of $\sqrt{T}(\bar{X}_T - \mu)$ highly dependent on $\mu$ even when the $X_i$ are normal with known unit variance. This difficulty can be resolved by replacing $\widehat{F}$, from which the bootstrap method resamples, by a more versatile resampling family $\widehat{F}_\mu$. Specifically, assuming $\mathrm{Var}(X_i) = 1$, the unknown distribution $G$ of $X_i - \mu$ can be estimated by the empirical distribution $\widehat{G}$ of $X_i - \bar{X}_T$. Let $\widehat{F}_\mu(\cdot) = \widehat{G}(\cdot - \mu)$ so that $\widehat{F}_\mu$ has mean $\mu$, and let $\widehat{u}_\alpha(\mu)$ be the $\alpha$-quantile of the sampling distribution of $\sqrt{T^*}(\bar{X}_{T^*}^* - \mu)$, where the $X_i^*$ are i.i.d. random variables from $\widehat{F}_\mu$. By analogy with the exact confidence set $\{\mu : u_\alpha(\mu) \leq \sqrt{T}(\bar{X}_T - \mu) \leq u_{1-\alpha}(\mu)\}$ proposed by Rosner and Tsiatis [17] when the $X_i$ are normal, for which the quantiles $u_\alpha(\mu)$ and $u_{1-\alpha}(\mu)$ of the sampling distribution of $\sqrt{T}(\bar{X}_T - \mu)$ can be determined by recursive numerical integration, Chuang and Lai [6] define the "hybrid" confidence set

$$(1.6) \qquad \{\mu : \widehat{u}_\alpha(\mu) < \sqrt{T}(\bar{X}_T - \mu) < \widehat{u}_{1-\alpha}(\mu)\}$$

as a hybrid of the exact and bootstrap methods. The hybrid confidence set is shown to be second-order accurate by making use of an Edgeworth expansion involving a $k$-variate normal distribution.

This derivation of second-order accuracy requires $k$ to be fixed and breaks down in the case of fully sequential (instead of group sequential) procedures, for which Edgeworth-type expansions are considerably more complicated and involve, besides the usual cumulants of $X_i$ (or $X_i^*$), certain fluctuation-theoretic quantities that are related to the random walk $\{S_n\}$ (or $\{S_n^*\}$), as shown by Woodroofe and Keener [29] and Lai and Wang [14]. On the other hand, (1.1) is usually satisfied by fully sequential truncated tests such as those studied by Woodroofe and his collaborators. Although this implies that $\sqrt{T}(\bar{X}_T - \mu)$ is an asymptotic pivot in these fully sequential settings when $T$ becomes infinite, the finite-sample behavior of $\sqrt{T}(\bar{X}_T - \mu)$ still varies substantially with $\mu$, as will be shown in the simulation studies in Section 3. In the context of exponential families, the asymptotic theory and the numerical studies of Woodroofe and his collaborators, however, suggest that suitable bias correction of $\sqrt{T}(\bar{X}_T - \mu)$ can improve its pivotal nature substantially, making the sampling distribution much less dependent on $\mu$. In Section 2 we first develop bias-corrected pivots in a nonparametric setting and then make use of them to construct nonparametric bootstrap confidence intervals. Numerical results are given in Section 3, where we also compare the different approaches to constructing confidence intervals following fully sequential tests. Further discussion of these results and some concluding remarks are given in Section 4.

## 2. Bias correction for a modified pivot and bootstrap confidence intervals under optional stopping

In this section we first review Woodroofe's pivot [25, 26] that corrects $\sqrt{T}(\bar{X}_T - \mu)$ for optional stopping in the case of a normal population with unknown mean $\mu$ and known variance 1. We then extend this method to nonparametric problems in which the normal mean $\mu$ is replaced by smooth functions of mean vectors of possibly non-normal populations. Note that the quantity of interest here is an approximate pivot for constructing confidence intervals, rather than point estimates of $\mu$ for which Siegmund [18], Emerson and Fleming [12] and Liu and Hall [15] have introduced bias-corrected or unbiased estimators following sequential tests.

Suppose the stopping rule is of the form

$$（2.1a） \qquad T = \min\{n_0(a),\ \max(t_a, n_1(a))\},$$

where $n_0(a) \sim a/\epsilon_0$ and $n_1(a) \sim a/\epsilon_1$, with $0 < \epsilon_0 < \epsilon_1$, and

$$（2.1b） \qquad t_a = \inf\{n \geq 1 : ng(S_n/n) \geq a\},$$

in which $g$ is continuously differentiable. A naive pivot is

$$（2.2） \qquad R_0(\mu) = \sqrt{T}(\bar{X}_T - \mu).$$

Whereas $\sqrt{n}(\bar{X}_n - \mu)$ has mean 0 and variance 1 for a fixed sample size $n$, optional stopping affects the first two moments of $\sqrt{T}(\bar{X}_T - \mu)$. First note that $a/t_a \xrightarrow{P} g(\mu)$ and therefore $a/T_a \xrightarrow{P} \kappa(\mu) := \max\{\epsilon_0, \min(g(\mu), \epsilon_1)\}$. As shown by Woodroofe [25],

$$（2.3） \qquad ER_0(\mu) \doteq a^{-1/2}[(d/d\mu)\kappa^{1/2}(\mu)] = (\kappa(\mu)/a)^{1/2} b(\mu),$$

where $b(\mu) = [(d/d\mu)\kappa^{1/2}(\mu)]/\kappa^{1/2}(\mu) = \dot{\kappa}(\mu)/\{2\kappa(\mu)\}$. This suggests the bias-corrected pivot

$$（2.4） \qquad R_1(\mu) = R_0(\mu) - T^{-1/2} b(\bar{X}_T)$$

as an improvement over (2.2). Moreover, Woodroofe [22] has shown that

$$（2.5） \qquad P(R_1(\mu) \leq x) \doteq \Phi(x) - (2a)^{-1} x\phi(x)[(d/d\mu)\kappa^{1/2}(\mu)]^2$$

in a very weak sense, i.e., the integral of the left-hand side of (2.5) with respect to $\xi(\mu)\,d\mu$ has an asymptotic expansion given by that of the right hand side (see (1.3)), where $\phi$ and $\Phi$ denote the standard normal density and distribution function, respectively. This in turn yields

$$（2.6） \qquad \begin{aligned} \int x^2 \, dP(R_1(\mu) \leq x) &\doteq \int x^2 \phi(x) \\ &\quad + (2a)^{-1} x^2(x^2-1)\phi(x)[(d/d\mu)\kappa^{1/2}(\mu)]^2 \, dx \\ &= 1 + a^{-1}[(d/d\mu)\kappa^{1/2}(\mu)]^2. \end{aligned}$$

Since $\{ER_1^2(\mu)\}^{1/2} \doteq 1 + (2a)^{-1}\kappa(\mu)b^2(\mu)$ by (2.6) and since $T^{-1} \doteq \kappa(\mu)/a$, these calculations led Woodroofe [25] to the approximate pivot (1.5).

We next remove the assumption of normality on the $X_i$ which we also extend to $d$-dimensional vectors. Instead of the mean vector $\mu$, we consider more general smooth functions $h$ of $\mu$ while the stopping time $T$ is still assumed to be of the form (2.1a,b).

Let $X, X_1, X_2, \ldots$ be i.i.d. $d \times 1$ random vectors with $EX = \mu$, $\text{Cov}(X) = V$ and $E\|X\|^r < \infty$ for some $r > 3$. Let $h : \mathbf{R}^d \to \mathbf{R}$ be twice continuously differentiable in some neighborhood of $\mu$. Consider a stopping rule $T$ of the form (2.1a,b), in which $g : \mathbf{R}^d \to \mathbf{R}$ is continuously differentiable in some neighborhood of $\mu$. Suppose $\epsilon_0 < g(\mu) < \epsilon_1$. Then application of the strong law of large numbers in conjunction with Taylor's theorem yields

$$\sqrt{T}\{h(\bar{X}_T) - h(\mu)\} \doteq \sqrt{T}(\nabla h(\mu))'(\bar{X}_T - \mu) + \sqrt{T}(\bar{X}_T - \mu)'\nabla^2 h(\mu)(\bar{X}_T - \mu)/2$$
$$(2.7) \doteq \frac{1}{\sqrt{a}} g^{1/2}(S_T/T)(S_T - \mu T)'\nabla h(\mu)$$
$$+ \frac{1}{2\sqrt{T}}\{T(\bar{X}_T - \mu)'\nabla^2 h(\mu)(\bar{X}_T - \mu)\},$$

in which the last approximate equality follows from $Tg(S_T/T) \doteq a$ (ignoring overshoot) so that $\sqrt{T} \doteq \sqrt{a}/g^{1/2}(S_T/T) \doteq \{a/g(\mu)\}^{1/2}$. By Wald's lemma, $E\{g^{1/2}(\mu) \times (S_T - \mu T)'\nabla h(\mu)\} = 0$. Moreover,

$$(2.8) \qquad g^{1/2}(S_T/T) - g^{1/2}(\mu) \doteq \{(\nabla g(\mu))'(S_T - T\mu)\}/\{2g^{1/2}(\mu)T\}.$$

By Anscombe's theorem [1], $\sqrt{T}(\bar{X}_T - \mu) = (S_T - \mu T)/\sqrt{T}$ has a limiting $N(0, V)$ distribution. Combining (2.7) with (2.8) and taking expectations, it can be shown by uniform integrability arguments that

$$(2.9) \begin{aligned} &E[\sqrt{T}\{h(\bar{X}_T) - h(\mu)\}] \\ &= \frac{(\nabla g(\mu))'V\nabla h(\mu)}{2(ag(\mu))^{1/2}} + \frac{1}{2}\left(\frac{g(\mu)}{a}\right)^{1/2} \text{tr}(\nabla^2 h(\mu)V) + o(a^{-1/2}). \end{aligned}$$

The second term on the right-hand side of (2.9) follows from $E(Z'AZ) = \text{tr}(AV)$ if $A$ is a nonrandom matrix and $Z$ is a random vector with $E(ZZ') = V$.

The difficult part of the proof of (2.9) lies in the technical arguments related to uniform integrability. For the case $h(x) = x$, Aras and Woodroofe [2] have provided such arguments to develop asymptotic expansions for the first four moments of $S_T/T$. Particularly relevant to our present problem are their Propositions 1, 2 and Section 5, which we can modify and refine to prove (2.9) when $\epsilon_0 < g(\mu) < \epsilon_1$. The details are omitted here. For $g(\mu) < \epsilon_0$ (or $g(\mu) > \epsilon_1$), stopping occurs at $n_0(a)$ (or $n_1(a)$) with probability approaching 1 and uniform integrability can again be used to show that

$$(2.10) \qquad E[\sqrt{T}\{h(\bar{X}_T) - h(\mu)\}] = \frac{1}{2}(n_i(a))^{-1/2}\text{tr}(\nabla^2 h(\mu)V) + o(a^{-1/2}),$$

with $i = 0$ or $1$ according as $g(\mu) < \epsilon_0$ or $g(\mu) > \epsilon_1$. Since $\kappa(\mu) = \max\{\epsilon_0, \min(g(\mu), \epsilon_1)\}$, $\nabla \kappa^{1/2}(\mu) = \frac{1}{2}\nabla g(\mu)/(g(\mu))^{1/2}$ if $\epsilon_0 < g(\mu) < \epsilon_1$, and $\nabla \kappa^{1/2}(\mu) = 0$ if $g(\mu) < \epsilon_0$ or $g(\mu) > \epsilon_1$. Recalling that $1/n_0(\mu) \sim \epsilon_0/a$ and $1/n_1(\mu) \sim \epsilon_1/a$, we can combine (2.9) and (2.10) into

$$(2.11) \qquad E[\sqrt{T}\{h(\bar{X}_T) - h(\mu)\}] = b(\mu, V)(\kappa(\mu)/a)^{1/2} + o(a^{-1/2}),$$

where

$$(2.12) \qquad b(\mu, V) = (\nabla \kappa^{1/2}(\mu))'V\nabla h(\mu)/\kappa^{1/2}(\mu) + \text{tr}(\nabla^2 h(\mu)V)/2.$$

For the special case $d = 1, h(\mu) = \mu$ and $V = 1$, $b(\mu, V) = [(d/d\mu)\kappa^{1/2}(\mu)]/\kappa^{1/2}(\mu)$, which agrees with Woodroofe's [25] approximation for $E\{\sqrt{T}(\bar{X}_T - \mu)\}$ derived from very weak asymptotic expansions for normal $X$.

Since $\mu$ is unknown, replacing $\mu$ by $\bar{X}_T$ and $\kappa(\mu)/a$ by $1/T$ in the last term of (2.3) leads to Woodroofe's [25] bias-corrected pivot $R_1(\mu)$ in (2.4). In the present nonparametric setting, $V$ is typically also unknown and has to be estimated to define both the naive pivot $R_0(\mu)$ and its bias-corrected version $R_1(\mu)$. Using the consistent estimates

$$(2.13) \quad \widehat{V}_T = \sum_{i=1}^T (X_i - \bar{X}_T)(X_i - \bar{X}_T)'/(T-1), \quad \widehat{\sigma}_T^2 = (\nabla h(\bar{X}_T))' \widehat{V}_T \nabla h(\bar{X}_T)$$

of $V$ and the asymptotic variance $\sigma^2 := (\nabla h(\mu))' V \nabla h(\mu)$ of $\sqrt{T}\{h(\bar{X}_T) - h(\mu)\}$, define

$$(2.14) \quad R_0(\mu) = \sqrt{T}\{h(\bar{X}_T) - h(\mu)\}/\widehat{\sigma}_T,$$

$$(2.15) \quad R_1(\mu) = [\sqrt{T}\{h(\bar{X}_T) - h(\mu)\} - T^{-1/2}b(\bar{X}_T, \widehat{V}_T)]/\widehat{\sigma}_T,$$

where $b(\mu, V)$ is defined in (2.12).

For the case of normal mean with known variance 1, Woodroofe [25] further refined $R_1(\mu)$ by scaling it with an asymptotic approximation to the standard deviation that he derived by very weak expansions; see (2.6) and (1.5). In the nonparametric setting with unknown covariance matrix $V$ considered here, such refinements are considerably much more complicated. In particular, better approximations to the asymptotic standard error than $\widehat{\sigma}_T$ (which is derived by linearizing $h$ around $\mu$) are needed. We therefore forgo such refinements and simply use the bias-corrected pivot $R_1(\mu)$ instead. The following example, which deals with the same testing problem as that considered in Woodroofe's [25] simulation study, shows that there is not much loss in the quality of the normal approximation to the distribution of $R_1(\mu)$ in comparison with that of $R(\mu)$.

**Example 1.** Suppose $X$ is normal with mean $\mu$ and known variance 1. Let

$$g(x) = (2\delta)^{-1}(\delta^2 + x^2)\mathbf{1}_{\{|x|\leq\delta\}} + |x|\mathbf{1}_{\{|x|>\delta\}},$$

which is symmetric and continuously differentiable, with $\dot{g} = x/\delta$ for $0 \leq x \leq \delta$ and $\dot{g} = 1$ for $x > \delta$. Take $a = 9$, $n_0(a) = 72$ and $n_1(a) = 1$ in (2.1a), as in Section 4 of [25]. Table 1 gives the $\alpha$-quantiles of $R(\mu)$, $R_0(\mu)$ and $R_1(\mu)$, respectively, computed from 10000 simulations, over different values of $\mu$ ranging from 0 to 1 as in Woodroofe's study. In this known variance setting, $R(\mu), R_0(\mu)$ and $R_1(\mu)$ are defined by (1.5), (2.2) and (2.4). Without assuming the variance to be known, we can use the version (2.15) for $R_1(\mu)$, which we denote by $R_{1,\widehat{\sigma}}(\mu)$ to indicate that the variance $V(=\sigma^2)$ is replaced by the sample variance $\widehat{\sigma}_T^2$. Table 1 also gives the $\alpha$-quantile of $R_{1,\widehat{\sigma}}(\mu)$ and the standard normal quantiles $z_\alpha$ for comparison. It shows that when $\alpha \leq 10\%$ or $\alpha \geq 90\%$, $q_\alpha(R)$ does not differ much from $z_\alpha$ and does not change much as $\mu$ varies between 0 and 1. The normal approximation is somewhat worse for $q_\alpha(R_1)$ or $q_\alpha(R_{1,\widehat{\sigma}})$ which, however, still does not change much as $\mu$ varies between 0 and 1. The normal approximation deteriorates substantially for $q_\alpha(R_0)$; moreover, $q_\alpha(R_0)$ is also markedly more variable with $\mu$, making it less "pivotal".

Since the stopping rule $T$ is assumed to be of the form (2.1a,b), $T/a \xrightarrow{P} 1/\kappa(\mu)$ and (1.1) is clearly satisfied. Therefore, by Anscombe's theorem, $R_0(\mu)$ is an asymptotic pivot and so is $R_1(\mu)$ that introduces a correction of the order $O_p(a^{-1/2})$ for

TABLE 1.

*Quantiles $z_\alpha, q_\alpha(R), q_\alpha(R_1), q_\alpha(R_{1,\widehat{\sigma}})$ and $q_\alpha(R_0)$ of the standard normal distribution, $R(\mu)$, $R_1(\mu)$, $R_{1,\widehat{\sigma}}(\mu)$ and $R_0(\mu)$, respectively*

| $\alpha$ (in %) | 2.5 | 5 | 10 | 20 | 50 | 80 | 90 | 95 | 97.5 |
|---|---|---|---|---|---|---|---|---|---|
| $z_\alpha$ | $-1.96$ | $-1.645$ | $-1.28$ | $-0.84$ | 0 | 0.84 | 1.28 | 1.645 | 1.96 |
| (a) $\mu = 0$ | | | | | | | | | |
| $q_\alpha(R)$ | $-1.9327$ | $-1.6280$ | $-1.2902$ | $-0.9053$ | 0.0035 | 0.8878 | 1.2832 | 1.6265 | 1.9880 |
| $q_\alpha(R_1)$ | $-1.9904$ | $-1.6581$ | $-1.3344$ | $-0.9372$ | $-0.0223$ | 0.9491 | 1.3206 | 1.6498 | 1.9721 |
| $q_\alpha(R_{1,\widehat{\sigma}})$ | $-2.0007$ | $-1.6929$ | $-1.3351$ | $-0.9163$ | 0.0045 | 0.9205 | 1.3012 | 1.7062 | 2.0083 |
| $q_\alpha(R_0)$ | $-2.2449$ | $-1.9255$ | $-1.6466$ | $-1.2218$ | 0.0075 | 1.1677 | 1.6419 | 1.9464 | 2.2369 |
| (b) $\mu = 0.25$ | | | | | | | | | |
| $q_\alpha(R)$ | $-1.9913$ | $-1.6558$ | $-1.2373$ | $-0.8129$ | 0.0029 | 0.8257 | 1.2795 | 1.6803 | 1.9698 |
| $q_\alpha(R_1)$ | $-2.0807$ | $-1.5996$ | $-1.2508$ | $-0.8376$ | $-0.0322$ | 0.8413 | 1.3063 | 1.6397 | 1.9515 |
| $q_\alpha(R_{1,\widehat{\sigma}})$ | $-1.9457$ | $-1.6106$ | $-1.2441$ | $-0.8274$ | $-0.0263$ | 0.8462 | 1.3106 | 1.6641 | 2.0033 |
| $q_\alpha(R_0)$ | $-1.9320$ | $-1.4949$ | $-1.0682$ | $-0.5817$ | 0.2862 | 1.0911 | 1.5139 | 1.8418 | 2.1377 |
| (c) $\mu = 0.5$ | | | | | | | | | |
| $q_\alpha(R)$ | $-1.8817$ | $-1.6496$ | $-1.2685$ | $-0.8450$ | $-0.0091$ | 0.8208 | 1.2448 | 1.6156 | 1.9321 |
| $q_\alpha(R_1)$ | $-2.0758$ | $-1.7250$ | $-1.3613$ | $-0.8840$ | $-0.0235$ | 0.8333 | 1.3260 | 1.6396 | 1.9818 |
| $q_\alpha(R_{1,\widehat{\sigma}})$ | $-2.0272$ | $-1.7000$ | $-1.3314$ | $-0.8745$ | $-0.0209$ | 0.8374 | 1.2769 | 1.6494 | 1.9985 |
| $q_\alpha(R_0)$ | $-1.7005$ | $-1.4107$ | $-1.0305$ | $-0.6237$ | 0.1894 | 1.0435 | 1.4668 | 1.8104 | 2.1436 |
| (d) $\mu = 0.75$ | | | | | | | | | |
| $q_\alpha(R)$ | $-1.9435$ | $-1.6190$ | $-1.3004$ | $-0.8758$ | $-0.0188$ | 0.8284 | 1.3224 | 1.6094 | 1.9519 |
| $q_\alpha(R_1)$ | $-2.0479$ | $-1.6768$ | $-1.3027$ | $-0.8684$ | $-0.0089$ | 0.8555 | 1.3153 | 1.6574 | 1.9600 |
| $q_\alpha(R_{1,\widehat{\sigma}})$ | $-2.0041$ | $-1.7406$ | $-1.3314$ | $-0.8787$ | $-0.0035$ | 0.8228 | 1.2873 | 1.6292 | 1.9441 |
| $q_\alpha(R_0)$ | $-1.7241$ | $-1.4270$ | $-1.1190$ | $-0.6657$ | 0.1415 | 1.0247 | 1.4202 | 1.8004 | 2.1070 |
| (e) $\mu = 1$ | | | | | | | | | |
| $q_\alpha(R)$ | $-1.9281$ | $-1.6617$ | $-1.2650$ | $-0.8719$ | $-0.0257$ | 0.7956 | 1.2505 | 1.6780 | 1.9335 |
| $q_\alpha(R_1)$ | $-2.0028$ | $-1.6963$ | $-1.3056$ | $-0.8719$ | $-0.0412$ | 0.8106 | 1.2514 | 1.6840 | 1.9542 |
| $q_\alpha(R_{1,\widehat{\sigma}})$ | $-2.0062$ | $-1.6821$ | $-1.3314$ | $-0.8344$ | $-0.0190$ | 0.8424 | 1.2520 | 1.6669 | 1.9506 |
| $q_\alpha(R_0)$ | $-1.8202$ | $-1.5173$ | $-1.1161$ | $-0.6909$ | 0.1150 | 0.9394 | 1.4222 | 1.8004 | 2.0165 |

$R_0(\mu)$. Although bootstrap confidence intervals based on $R_0(\mu)$ or $R_1(\mu)$ are therefore asymptotically valid, Example 1 suggests that $R_1(\mu)$ is more "pivotal" and may therefore provide substantial improvements over $R_0(\mu)$ for the finite-sample coverage errors. Simulation studies comparing both types of bootstrap confidence intervals are given in Section 3 to confirm this. A theoretical comparison would involve higher-order asymptotic expansions. While Woodroofe's very weak expansions are not applicable to the present nonparametric setting, the much more complicated Edgeworth-type expansions of Woodroofe and Keener [29] and Lai and Wang [14] can still be applied and will be presented elsewhere.

## 3. Numerical comparisons of various confidence intervals following sequential tests

Let $T$ be a stopping rule of the form (2.1a,b). Based on the sample $X_1, \ldots, X_T$ of random size $T$, the normal confidence interval

$$(3.1) \qquad (h(\bar{X}_T) - z_{1-\alpha}\widehat{\sigma}_T/\sqrt{T},\ h(\bar{X}_T) - z_\alpha\widehat{\sigma}_T/\sqrt{T})$$

simply uses the normal quantiles $z_\alpha$ and $z_{1-\alpha}$ to approximate the corresponding quantiles of $R_0(\mu)$, invoking Anscombe's theorem for its asymptotic justification. Similarly we can apply the normal approximation to the $\alpha$- and $(1-\alpha)$-quantiles of $R_1(\mu)$, leading to the interval

$$(3.2) \qquad \left( h(\bar{X}_T) - \frac{\{z_{1-\alpha}\widehat{\sigma}_T + T^{-1/2}b(\bar{X}_T, \widehat{V}_T)\}}{\sqrt{T}},\ h(\bar{X}_T) - \frac{\{z_\alpha\widehat{\sigma}_T + T^{-1/2}b(\bar{X}_T, \widehat{V}_T)\}}{\sqrt{T}} \right).$$

Instead of approximating the quantiles of $R_0(\mu)$ or $R_1(\mu)$ by normal quantiles, we can approximate them by the quantiles of the sequential bootstrap sample $(X_1^*, \ldots, X_{T^*}^*)$ described in the penultimate paragraph of Section 1, leading to the bootstrap confidence intervals based on $R_0(\mu)$ or $R_1(\mu)$.

The second paragraph of Section 1 has reviewed previous works on the exact method. As described more generally by Chuang and Lai [7, p. 2], the exact method involves (i) a family of distributions $F_\theta$ indexed by a real-valued parameter $\theta$ and (ii) a statistic $r(\theta; T, X_1, \ldots, X_T)$ for every given value of $\theta$, called a *root*. Let $u_\alpha(\theta)$ be the $\alpha$-quantile of $r(\theta; T, X_1, \ldots, X_T)$ under $F_\theta$. An exact equal-tailed confidence set for $\theta$ with coverage probability $1 - 2\alpha$ is

$$(3.3) \qquad \{\theta : u_\alpha(\theta) < r(\theta; T, X_1, \ldots, X_T) < u_{1-\alpha}(\theta)\}.$$

For the normal mean example (with $\theta = \mu$) considered by Siegmund [18], $r(\mu; T, X_1, \ldots, X_T) = (T, S_T)$, for which he introduced a total ordering to define the $p$-quantile $u_p(\theta)$. An obvious alternative choice is $r(\mu; T, X_1, \ldots, X_T) = \sqrt{T}(\bar{X}_T - \mu)$ that has been considered by Rosner and Tsiatis [17].

The exact method applies only when there are no nuisance parameters. In practice, however, not only do parametric models usually involve nuisance parameters, but one may also have difficulties in coming up with realistic parametric models. Without distributional assumptions on $X$, a $1 - 2\alpha$ level bootstrap confidence interval for a functional $\theta(F)$ of the distribution $F$ of $X$, based on the root $r(\theta; T, X_1, \ldots, X_T)$, is of the form

$$(3.4) \qquad \{\theta : u_\alpha^* < r(\theta; T, X_1, \ldots, X_T) < u_{1-\alpha}^*\},$$

where $u_\alpha^*$ is the $\alpha$-quantile of the distribution of $r(\widehat{\theta}; T^*, X_1^*, \ldots, X_{T^*}^*)$ in which $\widehat{\theta} = \theta(\widehat{F})$ and $(X_1^*, \ldots, X_{T^*}^*)$ is a bootstrap sample with random size $T^*$ drawn from the empirical distribution $\widehat{F}$ of $(X_1, \ldots, X_T)$. The bootstrap confidence interval (3.4) is tantamount to replacing $u_\alpha$ and $u_{1-\alpha}$ in (3.3) by $u_\alpha^*$ and $u_{1-\alpha}^*$ when $r(\theta; T, X_1, \ldots, X_T)$ is an approximate pivot, so that the quantile $u_p^*$ evaluated under $\widehat{F}$ can approximate the quantile under the true distribution $F$.

The hybrid method mentioned in Section 1 is based on reducing the nonparametric family $\mathcal{F}$ containing $F$ to another family $\widehat{F}_\theta$, where $\theta = \theta(F)$ is the unknown parameter of interest. It is particularly useful in situations where the sampling distribution of the root $r(\theta; T, X_1, \ldots, X_T)$ may depend on $\theta$ but is approximately constant over $\{F \in \mathcal{F} : \theta(F) = \theta\}$, as in group sequential clinical trials studied by Chuang and Lai [6] and in possibly nonstationary first-order autoregressive models considered in Section 5 of Chuang and Lai [7]. Applying the exact method to the family $\{\widehat{F}_\theta\}$ yields the hybrid confidence set

$$(3.5) \qquad \{\theta : \widehat{u}_\alpha(\theta) < r(\theta; T, X_1, \ldots, X_T) < \widehat{u}_{1-\alpha}(\theta)\}.$$

Note that hybrid resampling is a generalization of bootstrap resampling that uses the singleton $\{\widehat{F}\}$ as the resampling family.

In the following simulation studies we compare the hybrid confidence interval (3.5); the bootstrap confidence intervals $\text{Boot}(R_0)$ and $\text{Boot}(R_1)$ that use $R_0$ and $R_1$, respectively, as the root in (3.4); their direct normal approximation counterparts (3.1) and (3.2), denoted by $\text{Normal}(R_0)$ and $\text{Normal}(R_1)$, respectively; the exact confidence interval (3.3); and Woodroofe's [25] interval, denoted by $\text{Normal}(R)$, that uses the normal approximation to the renormalized pivot $R(\mu)$ in (1.5) derived under the parametric model. The quantiles in the bootstrap (or hybrid) confidence intervals are computed from 1000 samples drawn from $\widehat{F}$ (or $\widehat{F}_\theta$).

**Example 2.** Let $X_1, X_2, \ldots$ be i.i.d. $N(\mu, 1)$ random variables and let $T$ be the stopping rule of the form (2.1a,b) with $g(x) = x^2/2, a = 4.5, n_1(a) = 15$ and $n_0(a) = 75$. This corresponds to [4] repeated significance test (RST) of $H_0 : \mu = 0$ that stops sampling at $T = \inf\{n \geq 15 : |S_n| \geq 3\sqrt{n}\} \wedge 75$. Table 2 gives the coverage errors of the upper (U) and the lower (L) confidence bounds for $\mu$ following the RST, constructed by the various methods reviewed above with the nominal coverage error $\alpha = 5\%$. Each result is based on 10000 simulations. It shows that the hybrid confidence limits have coverage errors similar to those of the exact method (which should be 5%, with departures from 5% due to the Monte Carlo sampling variability). Woodroofe's method using normal approximation for the pivot $R$ also works well, except for a case ($\mu = 1.2$). While the normal confidence interval (3.1) using $R_0$ as the pivot has inaccurate coverage for certain values of $\mu$, the confidence interval (3.2), which uses $R_1$ as the pivot, and the bootstrap confidence intervals $\text{Boot}(R_0)$ and $\text{Boot}(R_1)$ show substantial improvement.

**Example 3.** The "exact" method for constructing the confidence interval (3.3) requires precise specification of a one-parameter family $F_\theta$, which we have assumed to be $N(\theta, 1)$ in Example 2. Woodroofe's pivot $R(\mu)$ is also derived under such parametric assumption. On the other hand, the pivots $R_0(\mu)$ and $R_1(\mu)$ can be derived nonparametrically. Suppose the underlying distribution $F$ is actually a mixture of $N(\mu, 1)$ and $\mu + (\text{Exp}(1) - 1)$, putting mixing probability 0.2 on $N(\mu, 1)$ and 0.8 on the other component of the mixture distribution, where Exp(1) denotes the exponential distribution with mean 1. Although the mixture distribution still has mean $\mu$ and variance 1, the skewness of the exponential component adversely affects the coverage errors of the "exact" confidence interval that assumes normality and those of the confidence intervals based on the normal approximation to $R(\mu), R_0(\mu)$ and $R_1(\mu)$, respectively, as shown in Table 3 that uses the same stopping rule as that of Example 2, $\alpha = 5\%$ and 10000 simulations to compute each coverage error. The bootstrap confidence interval $\text{Boot}(R_1)$ shows substantial improvement and the hybrid confidence interval performs even better.

TABLE 2.
Coverage errors (in %) for confidence limits of a normal mean $\mu$

| | Exact | | Hybrid | | Normal($R_0$) | | Boot($R_0$) | | Normal($R_1$) | | Boot($R_1$) | | Normal($R$) | |
|---|---|---|---|---|---|---|---|---|---|---|---|---|---|---|
| $\mu$ | L | U | L | U | L | U | L | U | L | U | L | U | L | U |
| 0.0 | 5.29 | 5.19 | 4.85 | 5.64 | 5.51 | 5.49 | 3.55 | 3.44 | 4.95 | 5.25 | 6.14 | 6.28 | 5.51 | 5.49 |
| 0.2 | 5.19 | 4.79 | 5.44 | 4.76 | 11.47 | 4.77 | 9.40 | 3.78 | 5.86 | 5.33 | 5.61 | 4.23 | 5.01 | 4.77 |
| 0.4 | 4.96 | 4.67 | 5.26 | 5.14 | 5.68 | 4.61 | 5.75 | 4.54 | 5.55 | 4.85 | 5.64 | 4.54 | 5.66 | 4.61 |
| 0.6 | 5.02 | 5.04 | 5.43 | 4.68 | 5.01 | 2.61 | 5.12 | 2.96 | 4.76 | 5.27 | 5.12 | 6.22 | 5.01 | 4.57 |
| 0.8 | 5.08 | 4.75 | 4.69 | 5.20 | 5.06 | 2.77 | 5.10 | 3.44 | 4.90 | 5.50 | 5.10 | 4.14 | 5.06 | 4.54 |
| 1.0 | 5.13 | 5.14 | 4.98 | 4.80 | 5.11 | 3.18 | 4.86 | 4.69 | 4.91 | 6.35 | 4.86 | 4.77 | 5.11 | 5.30 |
| 1.2 | 4.86 | 5.07 | 4.73 | 5.29 | 4.86 | 3.66 | 4.97 | 7.09 | 5.07 | 3.59 | 4.97 | 3.26 | 4.86 | 3.66 |
| 1.6 | 4.84 | 5.03 | 5.30 | 4.79 | 4.87 | 5.02 | 4.78 | 6.40 | 4.84 | 5.11 | 4.78 | 6.32 | 4.87 | 5.02 |

**Example 4.** The stopping rule $T$ in Examples 2 and 3 is associated with the RST when the variance is known. In the case of unknown variance, an obvious modification is

$$(3.6) \qquad T = \inf\{n \geq 15 : |S_n|/\widehat{\sigma}_n \geq 3\sqrt{n}\} \wedge 75,$$

where $\widehat{\sigma}_n^2 = n^{-1}\sum_{i=1}^n (X_i - \bar{X}_n)^2$. Note that $T$ is still of the form (2.1a,b), with $a = 4.5$ and $g : \mathbf{R} \times (0, \infty) \to [0, \infty)$ defined by

$$g(\eta, b) = \begin{cases} \eta^2/\{2(b - \eta^2)\} & \text{if } b \geq \eta^2, \\ 0 & \text{otherwise,} \end{cases}$$

since $ng(\sum_1^n X_i/n, \sum_1^n X_i^2/n) = S_n^2/(2n\widehat{\sigma}_n^2)$. Let $\mu = EX$ and $\mu_j = EX^j$ for $j \geq 1$. Then the covariance matrix of $(X, X^2)'$ has $\mu_2 - \mu^2$ and $\mu_4 - \mu_2^2$ as its diagonal elements, and its off-diagonal elements are both equal to $\mu_3 - \mu\mu_2$. Moreover, $(\nabla g)(\mu, \mu_2) = (\mu\mu_2, -\mu^2/2)'/(\mu_2 - \mu^2)^2$. To construct confidence intervals for $\mu$ in the case of unknown variance, the bias-corrected pivot $R_1(\mu)$ in (2.15) can be computed easily by setting $\nabla h(\mu, \mu_2) = (1, 0)'$ in (2.12), which corresponds to $h(\mu, \mu_2) = \mu$. Table 4 gives the coverage errors of the confidence limits for $\mu$, using $R_0$ and $R_1$, respectively, and also by bootstrapping $R_1$. The quantiles $t(R_0)$ and $t(R_1)$ of the $t$-distribution with $T$ degrees of freedom for the pivots $R_0$ and $R_1$ are used in lieu of normal quantiles, following Woodroofe and Coad [28, Section 4]. Again, $\alpha = 5\%$ and each result is based on 10000 simulations. It shows considerable improvement of $\text{Boot}(R_1)$ over the confidence limits based on $t$-distribution approximations for the approximate pivots $R_0$ and $R_1$.

## 4. Conclusion

In their discussion on pp. 33–36 of [7], Woodroofe and Weng have given some comparative studies of the hybrid method with the approach based on very weak expansions for constructing confidence intervals following group sequential tests. Their simulation studies have shown that the expansions "do not work very well for the repeated significance tests in Example 1 of the paper" and that "with large horizons, however, expansions work very well for triangular tests in which $g(x) = \delta + |x|$." They also raised the issue concerning "robustness with respect to the normality assumption." Their approach based on very weak expansions was subsequently extended by Morgan [16] to more general group sequential tests. This paper continues their investigation in other directions. First we consider fully sequential instead of group sequential tests. Secondly, we show the robustness of the bias correction, which Woodroofe derived by very weak expansions, by rederiving it without the parametric assumption. Thirdly, instead of applying directly normal approximation to the bias-corrected pivot, we use it as the root for constructing bootstrap confidence intervals. Simple bootstrapping indeed substantially reduces the computational cost of the hybrid method. Our conclusion from the comparative studies in Section 3 is that at the expense of greater computational cost, the hybrid method still provides the most reliable confidence intervals. Moreover, it is also much more versatile and can handle complex statistical models, as recently shown by Lai and Li [13] for the problem of valid confidence intervals for the regression parameter of a proportional hazards model following time-sequential clinical trials with censored survival data.

TABLE 3.
Coverage errors in (%) for confidence limits of the mean $\mu$ of a mixture of normal and exponential distributions

| $\mu$ | Exact | | Hybrid | | Normal($R_0$) | | Boot($R_0$) | | Normal($R_1$) | | Boot($R_1$) | | Normal($R$) | |
|---|---|---|---|---|---|---|---|---|---|---|---|---|---|---|
| | L | U | L | U | L | U | L | U | L | U | L | U | L | U |
| 0.0 | 5.69 | 4.61 | 4.67 | 5.27 | 5.98 | 4.84 | 3.45 | 3.83 | 5.93 | 3.86 | 4.64 | 6.13 | 5.98 | 4.84 |
| 0.2 | 7.43 | 4.52 | 4.95 | 6.17 | 13.70 | 4.49 | 10.42 | 4.60 | 8.65 | 4.62 | 6.51 | 4.61 | 7.39 | 4.49 |
| 0.4 | 7.15 | 4.59 | 4.84 | 6.37 | 7.96 | 4.49 | 6.40 | 5.01 | 7.69 | 4.07 | 6.30 | 5.01 | 7.63 | 4.49 |
| 0.6 | 6.15 | 4.74 | 4.82 | 6.09 | 6.13 | 2.49 | 5.52 | 3.12 | 6.53 | 4.56 | 5.52 | 6.47 | 6.13 | 4.29 |
| 0.8 | 6.32 | 4.84 | 5.09 | 5.73 | 6.32 | 2.58 | 5.53 | 3.69 | 6.33 | 4.92 | 5.53 | 4.63 | 6.32 | 4.62 |
| 1.0 | 5.88 | 4.06 | 5.10 | 5.72 | 5.86 | 2.34 | 5.16 | 3.63 | 6.20 | 4.79 | 5.16 | 4.61 | 5.86 | 4.22 |
| 1.2 | 6.63 | 3.93 | 5.18 | 5.95 | 5.63 | 2.33 | 4.97 | 5.38 | 6.17 | 3.16 | 4.97 | 2.37 | 5.63 | 2.33 |
| 1.6 | 5.83 | 3.90 | 5.06 | 5.85 | 5.85 | 3.90 | 5.29 | 4.66 | 5.90 | 3.10 | 5.29 | 4.66 | 5.85 | 3.90 |

TABLE 4.
Coverage errors (in %) of confidence limits for the mean $\mu$ of a normal distribution (left panel) and a mixture of normal and exponential distributions (right panel) when the variance is unknown

| | Normal distribution | | | | | | Mixed normal-exponential | | | | | |
|---|---|---|---|---|---|---|---|---|---|---|---|---|
| | $t(R_0)$ | | $t(R_1)$ | | Boot($R_1$) | | $t(R_0)$ | | $t(R_1)$ | | Boot($R_1$) | |
| $\mu$ | L | U | L | U | L | U | L | U | L | U | L | U |
| 0.0 | 5.85 | 5.87 | 6.02 | 6.12 | 4.98 | 3.66 | 3.41 | 12.10 | 3.98 | 11.56 | 3.51 | 6.71 |
| 0.2 | 12.89 | 5.25 | 7.56 | 5.19 | 5.54 | 3.52 | 7.26 | 8.02 | 4.00 | 7.42 | 5.81 | 3.26 |
| 0.4 | 5.59 | 5.04 | 6.37 | 4.80 | 6.42 | 4.82 | 3.64 | 6.94 | 3.02 | 6.87 | 5.98 | 4.83 |
| 0.6 | 5.32 | 3.03 | 5.88 | 5.75 | 5.75 | 5.48 | 2.52 | 6.51 | 2.58 | 9.10 | 5.64 | 6.45 |
| 0.8 | 5.80 | 3.59 | 5.83 | 6.55 | 5.90 | 5.52 | 2.45 | 9.77 | 2.69 | 10.41 | 6.50 | 6.32 |
| 1.0 | 5.69 | 4.58 | 5.52 | 6.40 | 4.89 | 5.13 | 2.43 | 10.66 | 2.50 | 10.91 | 5.25 | 5.62 |
| 1.2 | 5.49 | 5.28 | 5.80 | 5.85 | 5.65 | 4.20 | 2.46 | 10.36 | 2.82 | 10.67 | 5.54 | 5.38 |
| 1.6 | 5.46 | 5.64 | 5.55 | 5.74 | 6.11 | 4.69 | 2.34 | 10.33 | 2.34 | 10.70 | 5.38 | 4.96 |